\documentclass[preprint,12pt]{elsarticle}

\usepackage{fancyhdr}
\usepackage{amsthm}
\usepackage{amssymb}
\usepackage{amscd}
\usepackage{amsmath, pb-diagram}
\usepackage[all]{xy}
\usepackage{setspace}
\usepackage{graphics,epsfig}
\usepackage{graphicx}
\usepackage{float}
\usepackage{array}
\usepackage{xfrac}
\usepackage{faktor}
\usepackage{multicol}
\usepackage[mathscr]{euscript}
\usepackage{mathtools}
\usepackage{tikz-cd}
\usepackage[margin=2.9cm]{geometry}
\usepackage{url}
\usepackage{nicefrac}
\usepackage{xfrac}

\usepackage{enumerate}
\usepackage[hang,flushmargin]{footmisc}
\usepackage{amsmath}
\usepackage{amsfonts}
\usepackage{amssymb}
\usepackage{amsthm}
\usepackage{color}
\usepackage{blindtext}
\usepackage{tikz}
\usepackage{tikz-cd}
\usepackage{xcolor}
\usepackage{graphicx}
\usepackage{setspace}
\usepackage{float}
\usepackage{multicol}
\graphicspath{ {./images/} }
\usepackage{fullpage}
\usepackage{hyperref}
\usepackage{array}
\newcolumntype{C}[1]{>{\centering\arraybackslash}m{#1}}

\usepackage{amssymb}
\usepackage{amsthm}

\theoremstyle{plain}
\newtheorem{theo}{Theorem}[section]

\theoremstyle{definition}

\newtheorem{ex}[theo]{Example}

\theoremstyle{plain}
\newtheorem{theorem}{Theorem}[section]
\newtheorem{lemma}[theorem]{Lemma}

\newtheorem{proposition}[theorem]{Proposition}

\theoremstyle{definition}
\newtheorem{definition}[theorem]{Definition}

\newtheorem{?}[theorem]{Question}
\newtheorem{example}[theorem]{Example}

\journal{Discrete Mathematics}

\begin{document}

\begin{frontmatter}



\title{Generalized graph splines and the Universal Difference Property}


\author[Selma address]{Selma Alt{\i}nok}
\ead{sbhupal@hacettepe.edu.tr}
\affiliation[Selma address]{organization={Department of Mathematics, Hacettepe University }, %Department and Organization
            addressline={06800}, 
            city={Beytepe},
            state={Ankara},
            country={Turkey}}

\author[Katie address]{Katie Anders\corref{cor1}}
\cortext[cor1]{corresponding author}
\ead{kanders@uttyler.edu}
\affiliation[Katie address]{organization={Department of Mathematics, University of Texas at Tyler}, %Department and Organization
            addressline={3900 University Blvd.}, 
            city={Tyler},
            state={TX},
            postcode={75707}, 
            country={USA}}
            
\author[DAundergrad]{Daniel Arreola\fnref{Arreola grad address}}
\ead{darre@iastate.edu}
\fntext[Arreolad grad address]{Present address: Department of Mathematics, Iowa State University, 411 Morrill Rd., Ames, IA, 50011, USA}
\affiliation[DAundergrad]{organization={Department of Mathematics, California State University, Long Beach}, %Department and Organization
            addressline={1250 Bellflower Blvd.}, 
            city={Long Beach},
            state={CA},
            postcode={90840}, 
            country={USA}}

\author[Luisa address]{Luisa Asencio\fnref{Asencio grad address}}
\ead{lasencio@wustl.edu}
\fntext[Asencio grad address]{Present address: Department of Mechanical Engineering \& Materials Science, Washington University in St. Louis, 1 Brookings Dr., St. Louis, MO, 63130, USA}
\affiliation[Luisa address]{organization={Department of Mathematics, St. Edward's University}, %Department and Organization
            addressline={3001 South Congress}, 
            city={Austin},
            state={TX},
            postcode={78704}, 
            country={USA}}

\author[Chloe address]{Chloe Ireland\fnref{Ireland grad address}}
\ead{ireland.100@osu.edu}
\fntext[Ireland grad address]{Present address: Department of Mathematics, The Ohio State University, 231 West 18th Ave., Columbus, OH, 43210, USA}
\affiliation[Chloe address]{organization={Department of Mathematical Sciences, Carnegie Mellon University}, %Department and Organization
            addressline={Wean Hall 6113}, 
            city={Pittsburgh},
            state={PA},
            postcode={15213}, 
            country={USA}}
            
\author[Samet address]{Samet Sar{\i}o\u{g}lan}
\ead{ssarioglan@hacettepe.edu.tr}
\affiliation[Samet address]{organization={ Department of Mathematics, Hacettepe University}, %Department and Organization
            addressline={06800}, 
            city={Beytepe},
            state={Ankara},
            country={Turkey}}
            
\author[Luke address]{Luke Smith\fnref{Smith grad address}}
\ead{luke.smith@utdallas.edu}
\fntext[Smith grad address]{Present address: Department of Mathematical Sciences, University of Texas at Dallas, 800 W. Campbell Rd, Richardson, TX, 75080, USA}
\affiliation[Luke address]{organization={Department of Mathematics, University of Texas at Tyler}, %Department and Organization
            addressline={3900 University Blvd.}, 
            city={Tyler},
            state={TX},
            postcode={75707}, 
            country={USA}}

\begin{abstract}
We study the generalized graph splines introduced by Gilbert, Tymoczko, and Viel and focus on an attribute known as the Universal Difference Property (UDP). We prove that paths, trees, and cycles satisfy UDP.  We explore UDP on graphs pasted at a single vertex and use Prüfer domains to illustrate that not every edge labeled graph satisfies UDP.  We show that UDP must hold for any edge labeled graph over a ring $R$ if and only if $R$ is a Prüfer domain.  Lastly, we prove that UDP is preserved by isomorphisms of edge labeled graphs.
\end{abstract}




\end{frontmatter}


\section{Introduction}
\label{intro}
Splines are perhaps best known for their usage in analysis and for their applications in finding approximate solutions to differential equations, but splines also appear in a variety of other contexts including geometry and topology.  To unify these various notions of splines, Gilbert, Tymoczko, and Viel defined generalized splines on graphs in \cite{Gilbert}, and it is these generalized graph splines we consider here.

\begin{definition}
Given a graph $G= (V,E)$ and a commutative ring $R$, an \textit{edge labeling of $G$} is a function $\alpha: E \to \mathcal{I}(R)$, where $\mathcal{I}(R)$ is the set of all ideals in $R$. The pair $(G,\alpha)$ denotes a graph $G$ with edge labeling $\alpha$.
\end{definition}

\noindent
Because we will be considering edge labeled graphs, we will only work with commutative rings.

\begin{definition}
A \textit{generalized spline} on $(G,\alpha)$  over a ring $R$ is a function $\rho: V \to R$ such that for each edge $uv$, the difference $\rho(u) - \rho(v)$ is an element of $\alpha(uv)$. 
\end{definition}

Note that when working with a graph $G$ and a ring $R$, an edge labeling of $G$ labels the edges of $G$ with ideals in $R$, while a generalized spline labels the vertices of $G$ with elements of $R$.  Simply put, edges are labeled with ideals, while vertices are labeled with ring elements.  The ring $R$ is called the base ring for the edge labeled graph $(G, \alpha)$, and henceforth we shall refer to generalized splines simply as splines.  When working with an edge labeled graph $(G,\alpha)$ with $|V(G)|=n$, we often number the vertices of $G$ as $v_1, v_2, \ldots, v_n$ and denote a spline $\rho$ on $G$ by $\rho=\left(\rho(v_1), \rho(v_2),\ldots, \rho(v_n)\right)$.

\begin{ex}\label{spline on 3cycle}
Figure \ref{3cycle} illustrates an edge labeled graph with base ring $\mathbb{Z}$.  We see that $f=(2,16, 26)$ is a spline on the edge labeled $3$-cycle.  
\begin{figure}[H]
\begin{center}
\scalebox{0.2}{\includegraphics{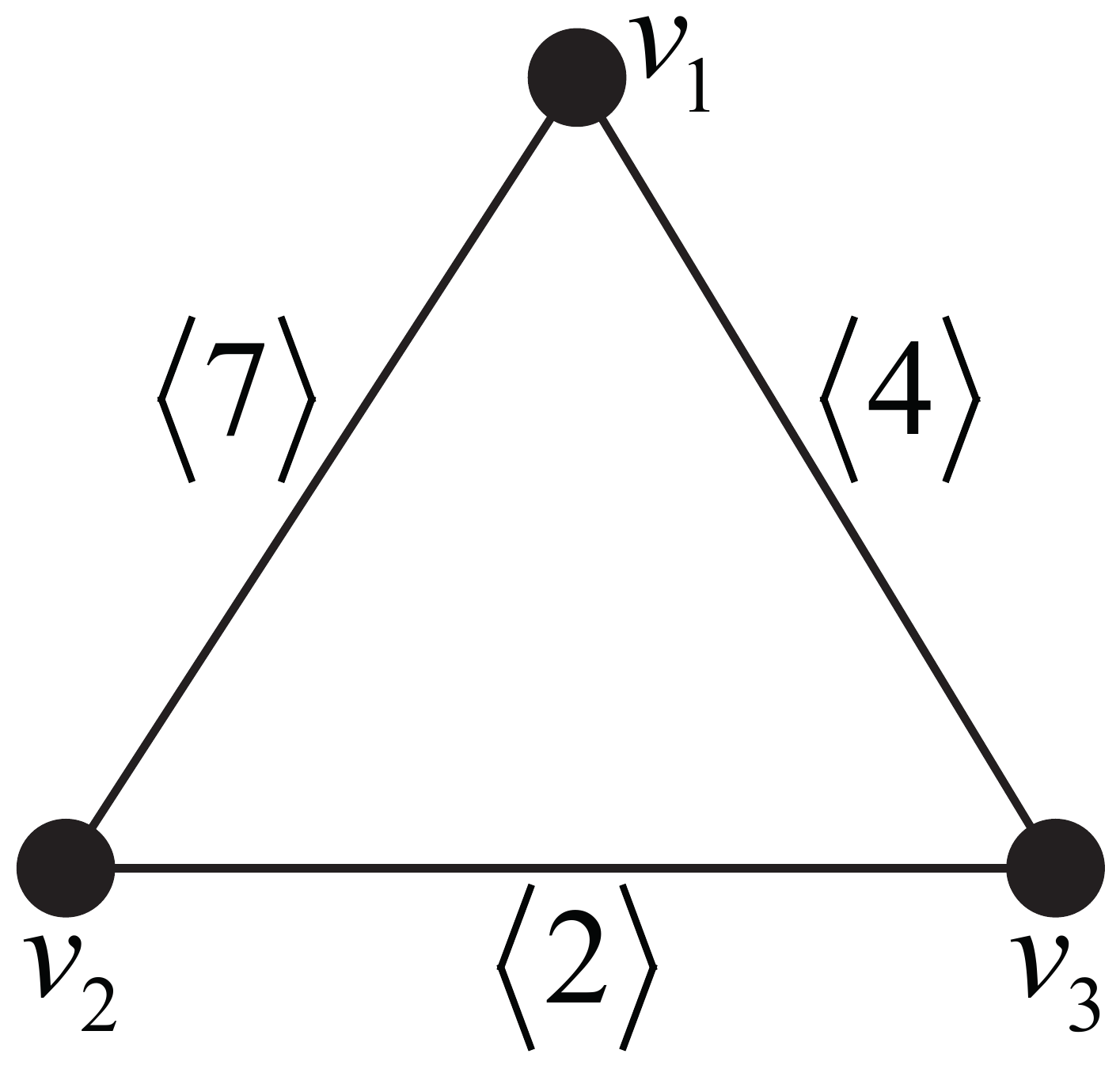}}
\caption{An edge labeled $3$-cycle}
\label{3cycle}
\end{center}
\end{figure}
\end{ex}

One of the primary results of \cite{Gilbert} is given in Proposition \ref{space of splines is a ring and an R-module}.

\begin{proposition}\label{space of splines is a ring and an R-module}
\cite{Gilbert} Fix a ring $R$ and a graph $G$ with edge labeling $\alpha$. The set of all splines on $(G,\alpha)$ is both a ring and an $R$-module.  We denote this set by $R_{(G, \alpha)}$.
\end{proposition}

The ring structure of $R_{(G,\alpha)}$ arises in algebraic and geometric topology. Equivariant cohomology rings of a GKM-manifold $X$ correspond to the generalized spline ring on the moment graph of $X$ over the multivariate polynomial ring with complex coefficients \cite{GKM}.  Given an edge labeled graph $(G,\alpha)$ over a ring $R$, a central question in spline theory is whether the set of splines $R_{(G, \alpha)}$ is free as an $R$-module.  When it is free, often the next task is to find a basis for $R_{(G,\alpha)}$.  The module structure of $R_{(G,\alpha)}$ and the freeness and bases of $R_{(G,\alpha)}$ have been studied by many mathematicians~\cite{Alt2,Alt,Bow,Tym, Dip,Gjo,Hand,Mah,Phi}.  Specifically, in addition to Proposition \ref{space of splines is a ring and an R-module}, another central result of \cite{Gilbert} is the description of a basis for the space of splines on an edge labeled tree.    Work has been done to find a basis for the space of splines on an edge labeled cycle, over $\mathbb{Z}$ in \cite{Bow} and \cite{Hand} and over $\mathbb{Z}/m\mathbb{Z}$ by Bowden and Tymoczko in \cite{Tym}.  In \cite{Rose and Suzuki}, Rose and Suzuki give a construction that can be used to find a basis for the basis of splines on any edge labeled graph over $\mathbb{Z}$.  Anders, Crans, Foster-Greenwood, Mellor, and Tymoczko investigated graphs admitting only constant splines in \cite{Anders}. It is worth noting that the freeness of $R_{(G,\alpha)}$ depends on the graph $G$ and the base ring $R$, but it does not depend on the ordering of the vertices of $G$.

Throughout this paper we use the following notation: for any $r \in R$, we let $\langle r \rangle$ denote the principal ideal of $R$ generated by $r$.  For connected vertices $u$ and $w$ in a graph $G$, we denote the set of all paths in $G$ from $u$ to $w$ by $\mathcal{P}_{(u,w)}$.

In \cite{Anders2}, Anders, Crans, Foster-Greenwood, Mellor, and Tymoczko proved the following theorem, noting that it enables us to compare the values of a spline on vertices that are in the same connected component of a graph.

\begin{theorem}\label{thm2.1}\cite{Anders2} Suppose $u, w \in V(G)$ for a graph $(G, \alpha)$, and $P = \langle u, v_1, \dots, v_n, w\rangle$ is a path from $u$ to $w$. Let $\alpha(P) = \alpha(uv_1) + \alpha(v_{1}v_2) + \dots + \alpha(v_nw)$. If $\rho \in R_{(G, \alpha)}$, then $\rho(u) - \rho(w) \in \alpha(P)$. Moreoever, if $P_1, P_2, \dots, P_m$ are paths from $u$ to $w$, then $\rho(u) - \rho(w) \in \bigcap_{i \in [m]} \alpha(P_i)$.
\end{theorem}

\begin{proof}
Consider a graph $(G, \alpha)$ with $u, w \in V(G)$, and let $P =  \langle u, v_1, \dots, v_n, w\rangle$ be a path from $u$ to $w$. Let $\alpha(P)$ be defined as above, and let $\rho$ be a spline on $(G, \alpha)$. Then
\[ \begin{split} \rho(u) - \rho(w) &= \rho(u) - \rho(v_1) + \rho(v_1) - \cdots - \rho(v_n) + \rho(v_n) - \rho(w) \\ &= (\rho(u) - \rho(v_1)) +( \rho(v_1) - \rho(v_2)) + \cdots + (\rho(v_{n-1}) - \rho(v_n)) + (\rho(v_n) - \rho(w)). \end{split} \]
Note that $\rho(u) - \rho(v_1) \in \alpha(uv_1),\; \rho(v_n) - \rho(w) \in \alpha(v_nw)$, and $\forall\, i \in [n-1]$, $\rho(v_i) - \rho(v_{i+1}) \in \alpha(v_iv_{i+1})$. Thus $\rho(u) - \rho(w) \in \alpha(P)$.

Now let $P_1, \dots, P_m$ be paths between $u$ and $w$. We know that $\forall\, i \in [m]$, $\rho(u) - \rho(w) \in \alpha(P_i)$, so $\rho(u) - \rho(w) \in \bigcap_{i \in [m]} \alpha(P_i)$. 
\end{proof}

The authors of \cite{Anders2} also asked when the converse to Theorem \ref{thm2.1} holds and posed the following question.

\begin{?}\label{when does UDP hold}
Let $(G, \alpha)$ be an edge labeled graph. Let $u, w$ be connected vertices in $(G, \alpha)$. Let $P_1, P_2, \dots, P_m$ be the paths in $G$ from $u$ to $w$. For each $i \in [m]$, let $\alpha(P_i)$ be the sum of the ideals labeling the edges in the path $P_i$. Let $x \in \bigcap_{i \in [m]} \alpha(P_i)$. Under what conditions does there exist a spline $\rho$ on $(G, \alpha)$ such that $\rho(u) - \rho(w) = x$?
\end{?}

To answer this question, the second, third, fourth, fifth, and seventh authors defined the Universal Difference Property, as described below.

\begin{definition}
If for every pair of vertices $u$ and $w$ and every $x \in \bigcap_{i \in [m]} \alpha(P_i)$ there exists a spline $\rho$ on $(G, \alpha)$ such that $\rho(u) - \rho(w) = x$, then we say that $G$ satisfies the \textit{Universal Difference Property}, abbreviated UDP.
\end{definition}

In Section \ref{paths trees cycles}, we show that UDP holds for edge labeled paths, trees, and cycles over any base ring $R$.  In Section \ref{pasting graphs and counterexample}, we consider taking two graphs on which UDP is satisfied and pasting them at a single vertex.  We determine a necessary and sufficient condition for UDP to hold on such a pasted graph.  We then turn our attention to Prüfer domains and find an example of an edge labeled, pasted graph over a Prüfer domain that does not satisfy UDP.  In Section \ref{UDP over a Prufer domain}, we prove that for a ring $R$, every edge labeled graph over $R$ satisfies UDP if and only if $R$ is a Prüfer domain.  We conclude in Section \ref{UDP and isomorphisms} by showing that the Universal Difference Property is a structural property of an edge labeled graph; that is, the Universal Difference Property is preserved by an isomorphism of edge labeled graphs.

\section{Paths, trees, and cycles}
\label{paths trees cycles}

In this section we answer Question \ref{when does UDP hold} for some basic types of graphs.  We begin with a lemma that will be used both in this section and in Section \ref{pasting graphs and counterexample}.

\begin{lemma}\label{LEMMA}
Consider an edge labeled graph $(G, \alpha)$.  Let $u, w \in V(G)$ and let $x \in \bigcap_{P \in \mathcal{P}_{(u, w)}} \alpha(P)$. Let $\varphi$ be a spline on $G$ such that $\varphi(u) - \varphi(w) = x$. Let $v \in V(G)$ and  $r \in R$. Then there exists a spline $\rho$ on $G$ such that $\rho(u) - \rho(w) = x$ and $\rho(v) = r$.
\end{lemma}

\begin{proof}
Let $u, w\in V(G)$. Let $x \in \bigcap_{P \in \mathcal{P}_{(u, w)}} \alpha(P)$. Then there exists a spline $\varphi$ on $G$ such that $\varphi(u) - \varphi(w) = x$. Let $v \in V(G)$ and $r \in R$. Consider the value $s = r - \varphi(v) \in R$. Create a new function $\rho$ on $V(G)$ such that $\forall \,z \in V(G),\, \rho(z) = \varphi(z) + s$. Note that $\rho$ is a spline on $G$ because $\varphi$ was. Furthermore, $\rho(u) - \rho(w) = \varphi(u) + s - \varphi(w)-s = x$. Finally, note that $\rho(v) = r$, as desired.
\end{proof}

This means that we always have one free variable when we are manually building splines. It is often convenient to set this free variable to be 0.  Now we will show that the Universal Difference Property must hold on any edge labeled path.

\begin{theorem} \label{theorem:3.1}
The Universal Difference Property holds when $G$ is a path.
\end{theorem}

\begin{proof}
Let $G$ be a path and let $u, w \in V(G)$. Because $G$ itself is a path, there is a unique path from $u$ to $w$. Denote it by $P_1 = \langle u, v_1, v_{2}, \ldots, v_{n}, w\rangle$. We have $\bigcap_{P \in \mathcal{P}_{(u, w)}} \alpha(P_i) = \alpha(P_1)$.  Let $x \in \alpha(P_1)$. Then by definition
$$x = \alpha_1 + \dots + \alpha_{n+1}, $$
where $\alpha_1 \in \alpha(uv_1),\, \alpha_{n+1} \in \alpha(v_nw),\, \text{ and }\alpha_{\ell} \in \alpha(v_{\ell - 1} v_{\ell})\text{ for all }2 \leq \ell \leq n$.  Let $r \in R$ be arbitrary and set $\rho(w) = r$. Now let $\rho(v_n) = \rho(w) + \alpha_{n+1}$. Note $\rho(v_n) - \rho(w) = \alpha_{n+1} \in \alpha(v_nw)$ by construction. For each $i \in [n-1]$, set $\rho(v_i) = \rho(w) + \sum_{j = i+1}^{n+1} \alpha_j$. Then for any $i \in [n-1]$, \[ \rho(v_i) - \rho(v_{i+1}) = \alpha_{i+1} \in \alpha(v_iv_{i+1}). \]
Finally, set $\rho(u) = \rho(w) + \sum_{j \in [n+1]} \alpha_j$. We have $\rho(u) - \rho(v_1) = \alpha_1 \in \alpha(uv_1)$. We can then find a spline for the rest of the path $G$ by giving the label $\rho(u)$ to all vertices coming before $u$ in the path $G$ and giving all of the vertices in $G$ after $w$ the label $\rho(w)$. Finally, note that 
\[ \rho(u) - \rho(w) = \rho(w) + \sum_{j \in [n+1]} \alpha_j - \rho(w) = \sum_{j \in [n+1]} \alpha_j = x. \]
We conclude that the Universal Difference Property holds when $G$ is a path.
\end{proof}

\begin{example}
Consider a path $P$ on $10$ vertices.  Suppose our ring is $\mathbb{Z}$ and the edges of $P$ are labeled as in Figure \ref{fig:path} . Note that there is only one path $P_1$ from $u$ to $w$. It has edge labels $\langle 4 \rangle, \langle 6 \rangle, \langle 2 \rangle, \langle 8 \rangle, \langle 12 \rangle,$ and $\langle 20 \rangle$. Thus $\alpha(P_1) = \langle 2 \rangle$. Suppose $x$ is arbitrarily chosen to be 64.  We can use the process described in the proof of Theorem \ref{theorem:3.1} to build a spline $\rho$ such that $\rho(u) - \rho(w) = x$.  Since we can write $64=8+12+4+8+12+20$, using the algorithm described in the proof of Theorem \ref{theorem:3.1}, we can label the vertices in the path from $u$ to $w$ as $64, 56, 44, 40, 32, 20, 0$, respectively.  All other vertices in $P$ are given the same label as their nearest neighbor in $P_1$.  Note that by construction we have $\rho(u)-\rho(w)=64-0=x$.

\begin{figure}[H]
\begin{center}
\scalebox{0.2}{\includegraphics{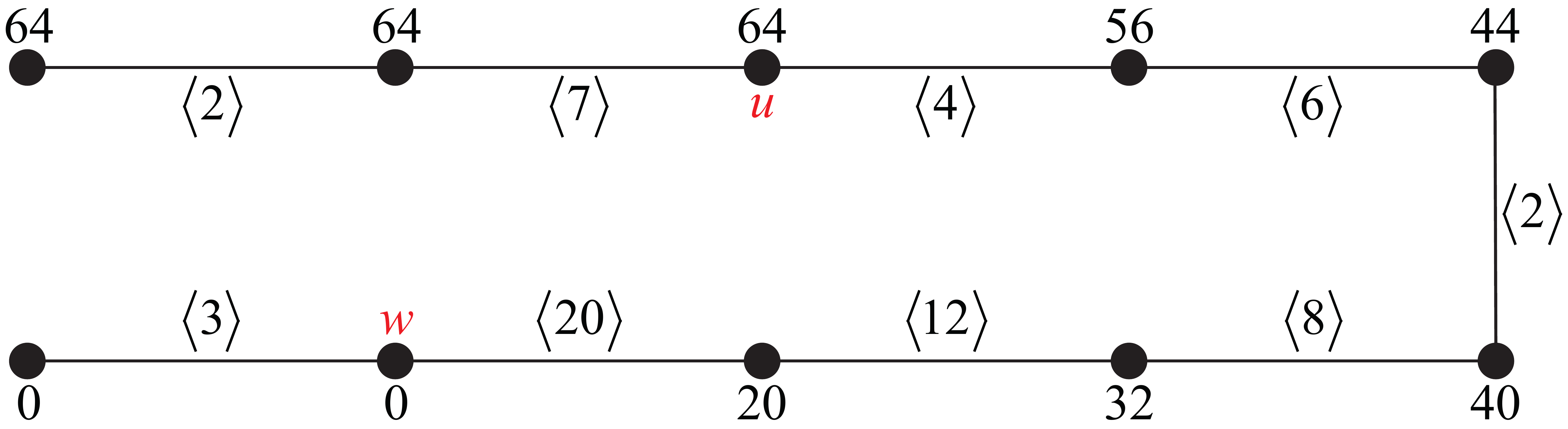}}
\caption{Path example}
\label{fig:path}
\end{center}
\end{figure}
\end{example}

\begin{theorem}
The Universal Difference Property holds when $G$ is a tree.
\end{theorem}

\begin{proof}
Let $G$ be a tree and $u, w\in V(G)$. Because $G$ is a tree, there exists a unique path $P_1$ between $u$ and $w$. Let $x \in \alpha(P_1)$ be arbitrary. One can create a spline $\rho$ on $G$ as follows.
As in Theorem \ref{theorem:3.1}, create a spline on $P_1$ such that $\rho(u) - \rho(w) = x$. Then let $\mathcal{U}$ be the set of vertices in $G$ that remain unlabeled. For any $v \in \mathcal{U}$, choose $z \in V(P_1)$ such that the unique path between $v$ and $z$ does not contain any elements of $P_1$ other than $z$. Assign $\rho(v) = \rho(z)$. We claim that this creates a spline. Indeed, let $a, b \in V(G)$ be adjacent. If $a, b \in V(P_1)$, by Theorem \ref{theorem:3.1} we know that $\rho(a) - \rho(b) \in \alpha(ab)$. If $a \in V(P_1)$ and $b \notin V(P_1)$, then $\rho(a) - \rho(b) = \rho(a) - \rho(a) = 0 \in \alpha(ab)$. The case where $a \notin V(P_1), b \in V(P_1)$ is similar. If $a, b \notin V(P_1)$, then $\rho(a) = \rho(b)$ because two adjacent vertices in $G$ not in $P_1$ must be closest to the same vertex on $P_1$. Otherwise, there would be a cycle in $G$. Then $\rho(a) - \rho(b) = 0 \in \alpha(ab)$. Because we have considered all cases, we may conclude that $\rho$ is a spline on the tree $G$.  Since $\rho(u) - \rho(w) = x$ by construction, we see that UDP holds on the tree $G$.  
\end{proof}

\begin{example}
Suppose our ring is $\mathbb{Z}$. Consider the labeled tree in Figure \ref{tree_example} below. The path $P_1$ from $u$ to $w$ has edge labels $\langle 5 \rangle, \langle 5 \rangle, \langle 11 \rangle$, and $\langle 12 \rangle$. Thus $\alpha(P_1) = \langle 1 \rangle$. Suppose $x$ is arbitrarily chosen to be $53$. We can write $53 = 15 + 15 + 11 + 12$, and so, using the algorithm described in the proof of Theorem \ref{theorem:3.1}, we can label the vertices in the path from $u$ to $w$ as $53, 38, 23, 12$, and $0$, respectively. All other vertices on the tree are then given the same label as their nearest neighbor in $P_1$. Notice by construction we have $\rho(u) - \rho(w) = 53 - 0 = x$.

\begin{figure}[H]
\begin{center}
\scalebox{0.2}{\includegraphics{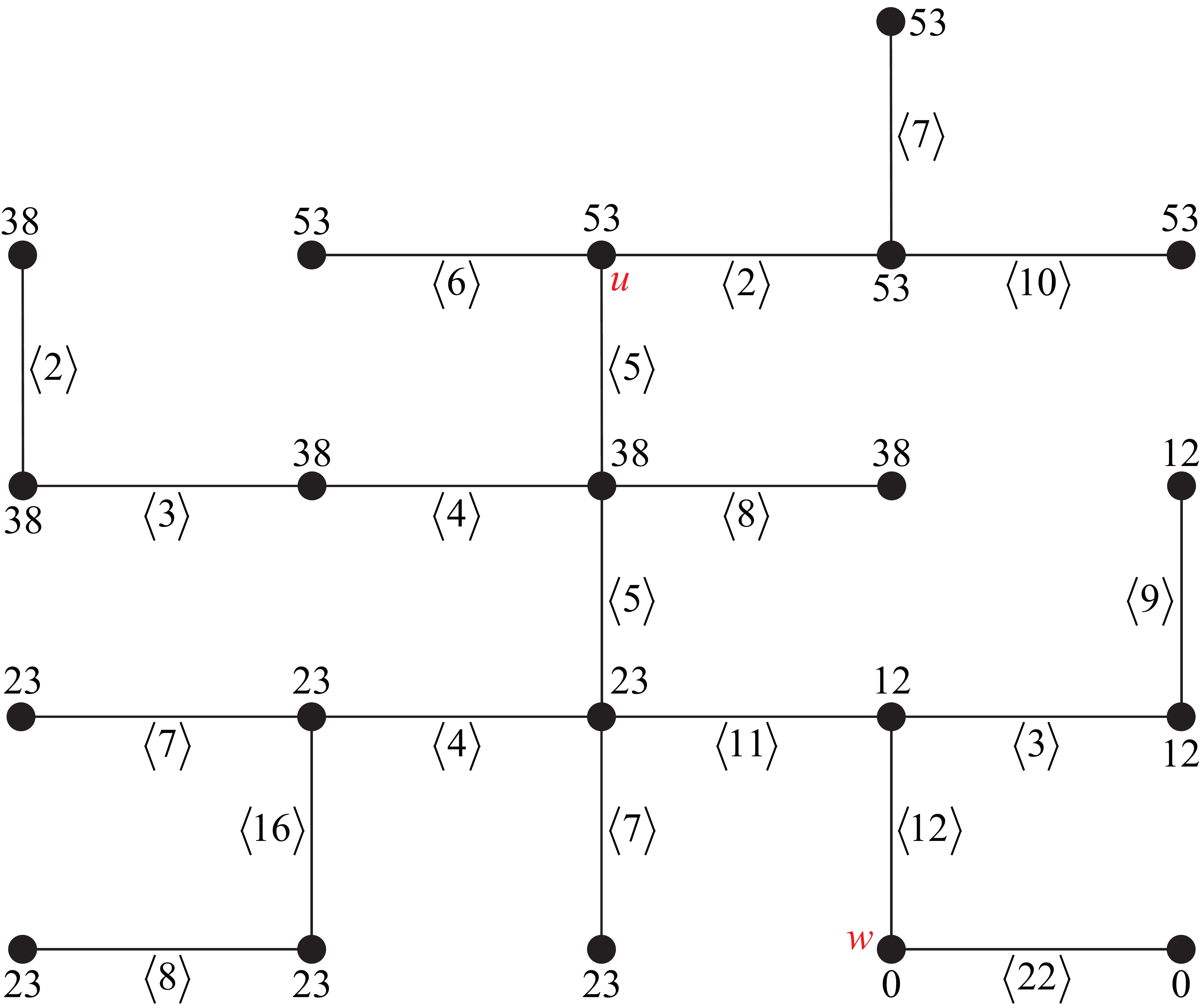}}
\caption{Tree example}
\label{tree_example}
\end{center}
\end{figure}
\end{example}

\begin{theorem}\label{UDP on cycles}
The Universal Difference Property holds when $G$ is a cycle.
\end{theorem}

\begin{proof}
Let $G$ be a cycle. There will be exactly two internally disjoint paths between any two distinct vertices $u$ and $w$. Denote them as
\[ P_1 = \langle u, v_1, \dots, v_n, w \rangle \]
and
\[ P_2 = \langle u, s_1, \dots, s_k, w \rangle. \]
Let $x \in \alpha(P_1) \cap \alpha(P_2)$. Then we know that
\begin{equation} x = \alpha_1 + \dots + \alpha_{n+1} \end{equation}
and
\begin{equation} x = \beta_1 + \dots + \beta_{k+1},  \end{equation}
for $\alpha_i, \beta_j$ in the appropriate ideals of paths $P_1, P_2$, respectively.

We construct a function $\rho: V \to R$ as follows. Let $r$ be an arbitrary element of the ring. Set $\rho(w) = r$ and let $\rho(v_n) =\rho(w)+ \alpha_{n+1}$. Note $\rho(v_n) - \rho(w) = \alpha_{n+1} \in \alpha(v_nw)$ by construction. For each $i \in [n-1]$, set $\rho(v_i) = \rho(w) + \sum_{j = i+1}^{n+1} \alpha_j$. Then for any $i \in [n-1]$,
\[ \rho(v_{i}) - \rho(v_{i+1}) = \left(\rho(w) + \sum_{j = i+1}^{n+1} \alpha_j\right)-\left(\rho(w) + \sum_{j = i+2}^{n+1} \alpha_j\right) = \alpha_{i+1}. \]
Finally, set $\rho(u) = \rho(w) + \sum_{j \in [n+1]} \alpha_j$, and observe that $\rho(u) - \rho(v_1) = \alpha_1 \in \alpha(uv_1)$ and $\rho(u) - \rho(w) = \sum_{i \in [n+1]} \alpha_i = x$.

Assign labels to the vertices $s_1, \dots, s_k$ by following the same process. However, now we have $\rho(u) = r + \sum_{i \in [n+1]} \alpha_i$ and $\rho(u) = r + \sum_{j \in [k+1]} \beta_j$. For this vertex labeling to be a function, these two values must be equal; indeed they are since the sums in these equations are (1) and (2), respectively, which both equal $x$. Thus $\rho$ is a spline on the cycle with $\rho(u)-\rho(w)=x$, and the Universal Difference Property holds for cycles.
\end{proof}

\begin{example}
Suppose our ring is $\mathbb{Z}$. Consider the edge labeling of the cycle in Figure \ref{fig:hendecagon} below. There are two paths between $u$ and $w$. The path $P_1$ has edge labels $\langle 2 \rangle, \langle 6 \rangle, \langle 8 \rangle, \langle 18 \rangle$, and $\langle 2 \rangle$. The path $P_2$ has edge labels $\langle 3 \rangle, \langle 5 \rangle, \langle 12 \rangle, \langle 11 \rangle, \langle 4 \rangle$, and $\langle 3 \rangle$. We have $\alpha(P_1) \cap \alpha(P_2) = \langle 2 \rangle + \langle 1 \rangle = \langle 1 \rangle$. Let $x = 48$. As described in the proof of Theorem \ref{UDP on cycles}, we can write $x$ in two different ways:
$x = 8 + 12 + 8 + 18 + 2$ and $x = 0 + 15 + 12 + 11 + 4 + 6$.
Labeling the vertices as shown in Figure \ref{fig:hendecagon} produces a spline $\rho$ on the cycle such that $\rho(u) - \rho(w) = 48 - 0 = x$.
\end{example}

%\begin{figure}[H]
%\centering
%\begin{tikzpicture}
%\draw (1.5,0) node{$\bullet$} node[below]{6} -- (3,0) node{$\bullet$} node[below]{6} %node[midway,above]{$\langle3\rangle$} -- (4.5,0) node{$\bullet$} node[below]{0} node[above]{$w$} %node[midway,above]{$\langle2\rangle$}
%-- (6,1.5) node{$\bullet$} node[right]{2} node[midway,left]{$\langle2\rangle$}
%-- (6,3) node{$\bullet$} node[right]{20} node[midway,left]{$\langle18\rangle$} -- (6,4.5) node{$\bullet$} %node[right]{28} node[midway,left]{$\langle8\rangle$}-- (4.5,6) node{$\bullet$} node[above]{40} %node[midway,left]{$\langle6\rangle$}
%-- (3,6) node{$\bullet$} node[above]{48} node[midway,below]{$\langle2\rangle$} node[below]{$u$} -- (1.5,6) %node{$\bullet$} node[above]{48} node[midway,below]{$\langle3\rangle$}
%-- (0,4.5) node{$\bullet$} node[left]{33} node[midway,right]{$\langle5\rangle$} -- (0,3)  node{$\bullet$} %node[left]{21} node[midway,right]{$\langle12\rangle$} -- (0,1.5) node{$\bullet$} node[left]{10} %node[midway,right]{$\langle11\rangle$} 
%-- (1.5,0) node[midway,right]{$\langle4\rangle$}; 
%\end{tikzpicture}
%\caption{Example 3.6}
%\label{fig:cycle_example}
%\end{figure}

\begin{figure}[H]
\begin{center}
\scalebox{0.23}{\includegraphics{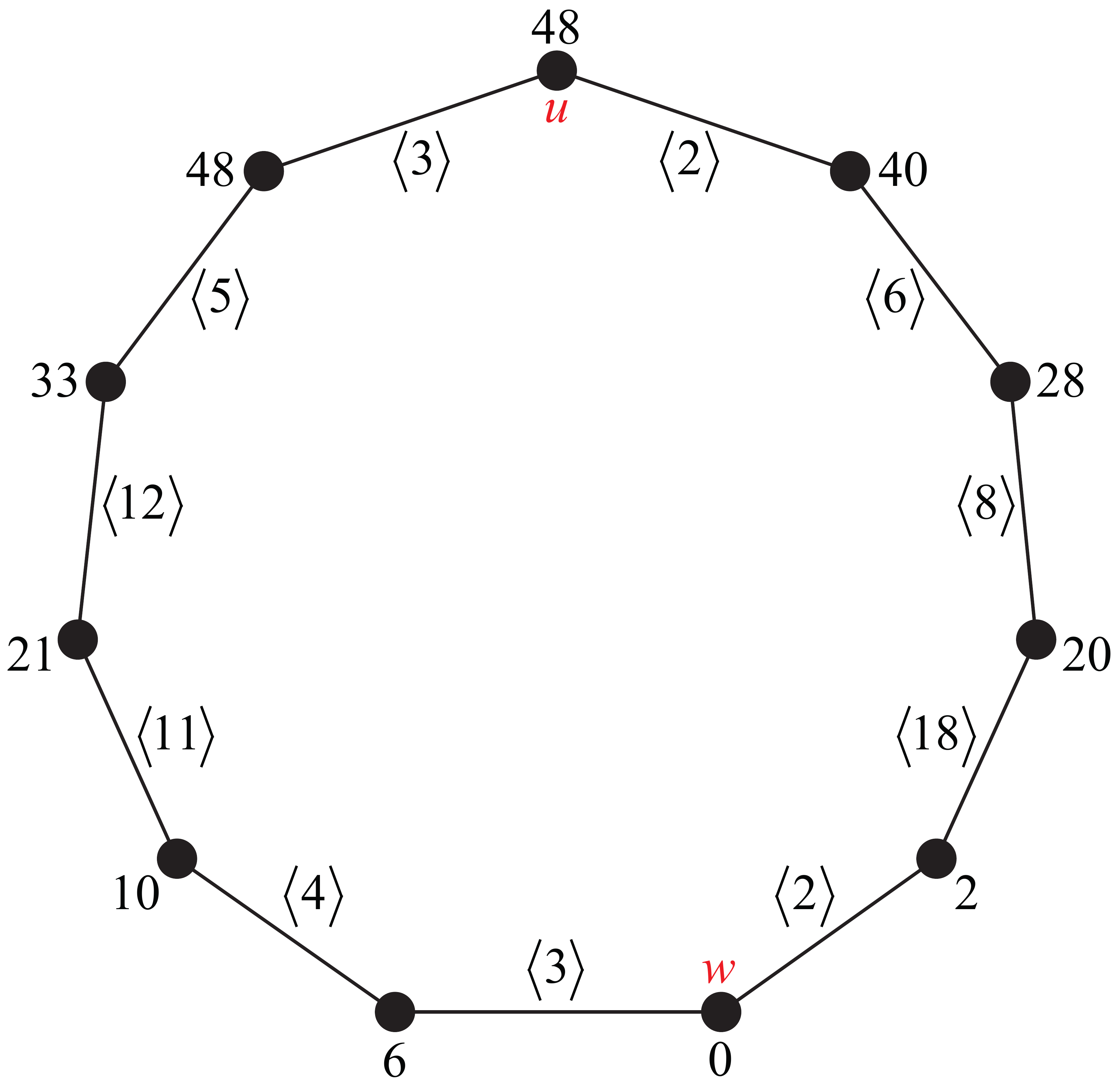}}
\caption{Cycle example}
\label{fig:hendecagon}
\end{center}
\end{figure}

\section{Pasting graphs and an example where UDP does not hold}
\label{pasting graphs and counterexample}
In this section, we examine when the Universal Difference Property holds on graphs pasted at a single vertex.  Specifically, we develop a necessary and sufficient condition for UDP to hold on an edge labeled graph constructed by pasting together at a single vertex two edge labeled graphs on which UDP holds.  Finally, we introduce Prüfer domains to allow us to provide an example of a graph on which UDP does not hold.  Our first step toward the necessary and sufficient condition is the following lemma.

\begin{lemma} \label{lemma:5.1} Let $G_1$ and $G_2$ be connected, edge labeled graphs satisfying the Universal Difference Property. Suppose $V(G_1) \cap V(G_2) = \{ z\}$. Let $G = G_1 \cup G_2$, where $G$ is created by pasting $G_1$ and $G_2$ at $z$. If $u$ and $w$ are both in $V(G_1)$ or both in $V(G_2)$, then for all $x \in \bigcap_{P \in \mathcal{P}_{(u, w)}} \alpha(P)$, there exists a spline $\rho$ on $G$ such that $\rho(u) - \rho(w) = x$. \end{lemma}
\begin{proof} 
Assume without loss of generality that $u, w \in V(G_1)$. Observe that no path from $u$ to $w$ can include an element of $V(G_2)$ other than $z$. If a path from $u$ to $w$ were to leave $G_1$, it would have to pass through $z$ by construction; but then to reach $w$, the path would have to return to $z$, violating the definition of a path. Thus each path from $u$ to $w$ involves only vertices in $G_1$. Since $G_1$ satisfies the Universal Difference Property, for any $x \in \bigcap_{P \in \mathcal{P}_{(u,w)}} \alpha(P)$ there exists a spline $\rho$ on $G_1$ such that $\rho(u) - \rho(w) = x$. Extend $\rho$ to the remaining vertices of $G$ by labeling all vertices in $G_2$ with $\rho(z)$. Then $\rho$ is a spline on $G$ satisfying $\rho(u) - \rho(w) = x$. \end{proof}

Lemma \ref{lemma:5.1} addresses the case where $u$ and $w$ both belong to $V(G_1)$ or both belong to $V(G_2)$.  What if exactly one of $u$ and $w$ is in $V(G_1)$ and the other is in $V(G_2)$?  The next theorem addresses this case, providing two conditions, and if either condition is satisfied, then UDP holds on $G_1\cup G_2$.

\begin{theorem}\label{theorem:5.2} Let $G_1$ and $G_2$ be connected, edge labeled graphs satisfying the Universal Difference Property. Suppose $V(G_1) \cap V(G_2) = \{ z\}$. Let $G = G_1 \cup G_2$, where $G$ is created by pasting $G_1$ and $G_2$ at $z$.  If, for every $u,w \in V(G_1 \cup G_2)$, either
\begin{equation}\label{uw in uz}
\bigcap_{P \in \mathcal{P}_{(u,w)}} \alpha (P) \subseteq \bigcap_{P \in \mathcal{P}_{(u,z)}} \alpha (P)
\end{equation}
or
\begin{equation}\label{uw in wz}
\bigcap_{P \in \mathcal{P}_{(u,w)}} \alpha (P) \subseteq \bigcap_{P \in \mathcal{P}_{(w,z)}} \alpha (P),
\end{equation}
then $(G_1 \cup G_2,\alpha)$ satisfies the Universal Difference Property.
\end{theorem}

\begin{proof}
Let $u, w \in V(G)$. The case where $u$ and $w$ are both in $V(G_1)$ or $u$ and $w$ are both in $V(G_2)$ is covered in Lemma \ref{lemma:5.1}. It remains to consider the case when $u$ and $w$ are on different connected graphs $G_1$ and $G_2$. Assume without loss of generality that $u \in V(G_1)$ and $w \in V(G_2)$ with $u$ and $w$ both distinct from $z$. Suppose $(3)$ holds and let $x \in \bigcap_{P \in \mathcal{P}_{(u,w)}} \alpha(P)$. Then $x \in \bigcap_{P \in \mathcal{P}_{(u,z)}} \alpha(P)$ by assumption. Since $G_1$ satisfies the Universal Difference Property, there exists a spline $\rho$ on $(G_1, \alpha|_{E(G_1)})$ such that $\rho(u) - \rho(z) = x$. This is important to note since every path from $u$ to $w$ in this case must contain $z$. Now extend the domain of $\rho$ to $G_1 \cup G_2$ by mapping every element in $V(G_2)$ to $\rho(z)$. It is clear that $\rho$ is a spline on $G_1 \cup G_2$ and since $w \in V(G_2)$, we have $\rho(u) - \rho(w) = \rho(u) - \rho(z) = x$ as desired. A similar argument proves the case when $(4)$ is true.
\end{proof}

\begin{example}
Suppose our ring is $\mathbb{Z}$. Consider the graph with edge labeling shown below in Figure \ref{pastingg}.  Note that $\bigcap_{P \in \mathcal{P}_{(u,z)}} \alpha(P) = \langle 2 \rangle$ and $\bigcap_{P \in \mathcal{P}_{(u,w)}} \alpha(P) = \langle 2 \rangle$. If we let $x = 38$, then we can then use the process described in Theorem \ref{theorem:5.2} to create a spline satisfying $\rho(u) - \rho(w) = 38 - 0 = x$.

\begin{figure}[H]
\begin{center}
\scalebox{0.2}{\includegraphics{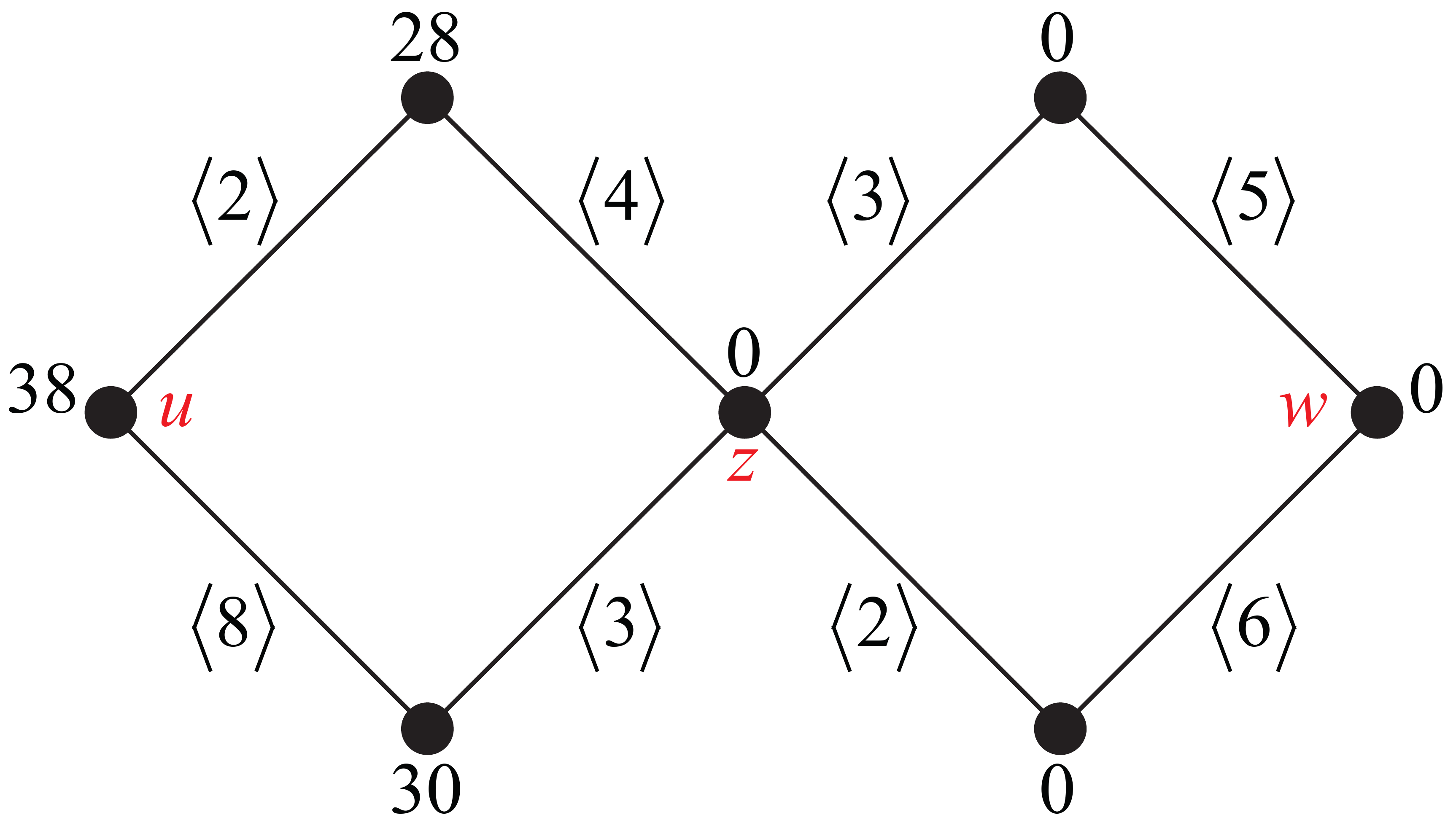}}
\caption{Example of an edge labeled graph created by pasting at a single vertex two edge labeled subgraphs each satisfying UDP}
\label{pastingg}
\end{center}
\end{figure}
\end{example}

The next lemma gives set containments that must hold for any graph $G=G_1\cup G_2$ pasted at a single vertex, regardless of whether $G_1$ and $G_2$ satisfy UDP.  In the proof of this next lemma, we see that the containments that are the reverses of \eqref{uw in uz} and \eqref{uw in wz} must always hold.

\begin{lemma} \label{prop:subset}
Suppose $V(G_1) \cap V(G_2) = \{ z\}$. Let $G = G_1 \cup G_2$, where $G$ is created by pasting $G_1$ and $G_2$ at $z$. Suppose $u \in V(G_1)$ and $w \in V(G_2)$. Then
$$\left(\bigcap_{P \in \mathcal{P}_{(u,z)}} \alpha(P)\right)  + \left(\bigcap_{P \in \mathcal{P}_{(z,w)}} \alpha(P)\right)  \subseteq \left(\bigcap_{P \in \mathcal{P}_{(u,w)}} \alpha(P)\right).$$
\end{lemma}

\begin{proof}
Let $\mathcal{P}_{(u,w)} = \{P_1, P_2, \ldots, P_m\}$. For each $P_i$, let $P_i'$ denote the subpath of $P_i$ from $u$ to $z$. Since $\alpha(P_i')$ is just a truncated sum of ideals in $\alpha(P_i)$ we have $\alpha(P_i') \subseteq \alpha(P_i)$ for every $1 \leq i \leq m$. It follows that $\bigcap_{i=1}^m \alpha(P_i') \subseteq \bigcap_{i=1}^m \alpha(P_i)$, or equivalently,
$$\bigcap_{P \in \mathcal{P}_{(u,z)}} \alpha(P) \subseteq \bigcap_{P \in \mathcal{P}_{(u,w)}} \alpha(P).$$
By symmetry we also have
$$\bigcap_{P \in \mathcal{P}_{(z,w)}} \alpha(P) \subseteq \bigcap_{P \in \mathcal{P}_{(u,w)}} \alpha(P).$$
Thus $$\left(\bigcap_{P \in \mathcal{P}_{(u,z)}} \alpha(P) \right) + \left(\bigcap_{P \in \mathcal{P}_{(z,w)}} \alpha(P) \right) \subseteq \left(\bigcap_{P \in \mathcal{P}_{(u,w)}} \alpha(P) \right).$$

\end{proof}

We are now ready to state and prove our main result on pasted graphs, Theorem \ref{pasting}, which gives a necessary and sufficient condition for UDP to hold on an edge labeled graph constructed by pasting together at a single vertex two edge labeled graphs on which UDP holds.  Theorem \ref{theorem:5.2} gives way to Theorem \ref{pasting}. Indeed, if the hypotheses of Theorem \ref{theorem:5.2} hold, then Equation \eqref{importantequation} is satisfied.

\begin{theorem}\label{pasting}
Let $G_1$ and $G_2$ be connected graphs such that the Universal Difference Property holds for each and $V(G_1) \cap V(G_2) = \{ z\}$. Let $G = G_1 \cup G_2$, where $G$ is created by pasting $G_1$ and $G_2$ at $z$. Then $G$ satisfies the Universal Difference Property if and only if for all $u \in V(G_1), w \in V(G_2)$, we have that 
\begin{equation}
\left(\bigcap_{P \in \mathcal{P}_{(u,w)}} \alpha(P) \right) = \left(\bigcap_{P \in \mathcal{P}_{(u,z)}} \alpha(P) \right) + \left(\bigcap_{P \in \mathcal{P}_{(z,w)}} \alpha(P) \right).
\label{importantequation} \end{equation}
\end{theorem}

\begin{proof}
Let $u, w \in V(G)$. The case where $u, w \in V(G_1)$ or $u, w \in V(G_2)$ is covered in Lemma \ref{lemma:5.1}. Hence, it remains to consider cases where, without loss of generality, $u \in V(G_1)$ and $w \in V(G_2)$. 

We will prove the forward direction by considering the contrapositive. Assume that Equation \eqref{importantequation} does not hold. By Lemma \ref{prop:subset}, we know it must be that the left side of (\ref{importantequation}) is not a subset of the right side. Then consider $x \in \bigcap_{P \in \mathcal{P}_{(u,w)}} \alpha(P)$ such that for all $s \in \bigcap_{P \in \mathcal{P}_{(u,z)}} \alpha(P)$ and for all $t \in\bigcap_{P \in \mathcal{P}_{(z, w)}} \alpha(P)$, $x \neq s + t$. Suppose for the sake of contradiction that there is a spline $\rho$ on $G$ such that $\rho(u) - \rho(w) = x$. By Theorem \ref{thm2.1}, we know that $\rho(u) - \rho(z) \in \left(\bigcap_{P \in \mathcal{P}_{(u,z)}} \alpha(P) \right)$; similarly, $\rho(z) - \rho(w)\in \left(\bigcap_{P \in \mathcal{P}_{(z,w)}} \alpha(P) \right)$. Let $s=\rho(u)-\rho(z)$ and $t=\rho(z)-\rho(w)$.  Then
\begin{align*} x &= \rho(u) - \rho(w) \\ &= \rho(u) - \rho(z) + \rho(z) - \rho(w) \\ &= s + t. \end{align*}
This yields a contradiction, since we chose $x$ such that it would be impossible to write it in this form. Hence our assumption was incorrect; there is no spline on $G$ such that $\rho(u) - \rho(w) = x$. Therefore the Universal Difference Property does not hold in this case.
 
 Assume that Equation \eqref{importantequation} holds. Let $x \in\bigcap_{P \in \mathcal{P}_{(u,w)}} \alpha(P)$. Then there exists $s \in\bigcap_{P \in \mathcal{P}_{(u,z)}} \alpha(P)$ and $t \in \bigcap_{P \in \mathcal{P}_{(z,w)}} \alpha(P)$ such that $x = s + t$. Because $s \in\bigcap_{P \in \mathcal{P}_{(u,z)}} \alpha(P)$, and $G_1$ satisfies the Universal Difference Property, there is a spline $\rho_{1}$ on $G_1$ such that $\rho_{1}(u) - \rho_{1}(z)= s$. Similarly, there is a spline $\rho_2$ on $G_2$ such that $\rho_2(z) - \rho_2(w) = t$, and by Lemma \ref{LEMMA} there is such a spline $\rho_2$ also satisfying $\rho_2(z) = \rho_1(z)$ since $\rho_1(z) \in R$. Now define the function $\rho$ on all of $G$ by \[ \rho(v) = \begin{cases} \rho_{1}(v), & v \in G_1 \\ \rho_{2}(v), & v \in G_2\text{.} \end{cases}\] Note that this indeed is a spline. Furthermore,
\[ \rho(u) - \rho(w) = \rho_{1}(u) - \rho_{1}(z) + \rho_{2}(z) - \rho_{2}(w) = s + t = x. \]
Therefore, combining this result and the result from Lemma \ref{lemma:5.1}, we know that for any $u, w\in V(G)$, for all $x \in \bigcap_{P \in \mathcal{P}_{(u, w)}} \alpha(P)$, there is a spline $\rho$ such that $\rho(u) - \rho(w) = x$ provided that (\ref{importantequation}) is true. Hence UDP is satisfied precisely when Equation \eqref{importantequation} holds. 
\end{proof}

\begin{example}
Below is an example of a pasted graph $(G,\alpha)$ satisfying Equation \eqref{importantequation}.

\begin{figure}[H]
\begin{center}
\scalebox{0.2}{\includegraphics{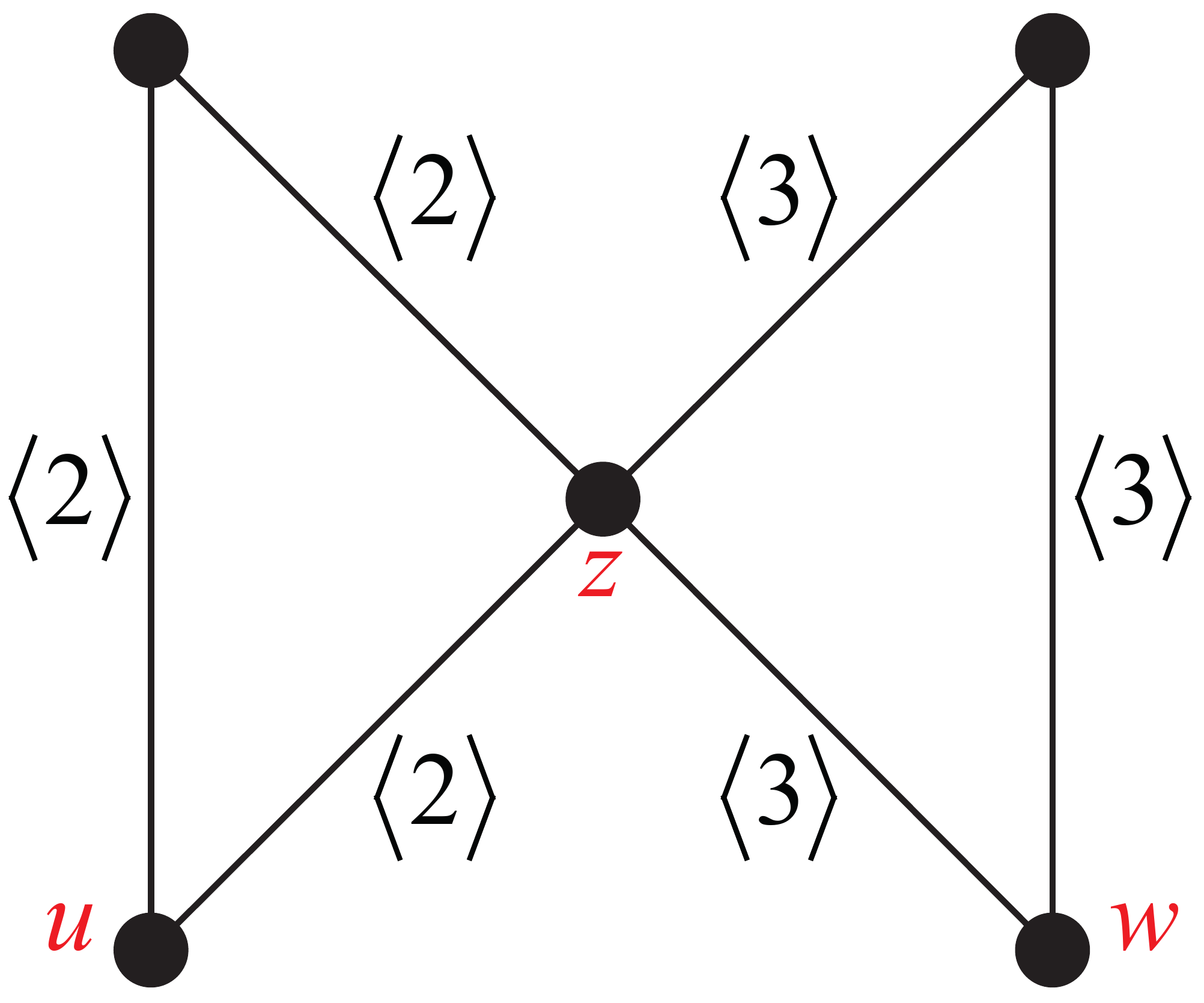}}
\caption{Example of a pasted graph that satisfies Equation \eqref{importantequation}}
\label{fig:bowtie}
\end{center}
\end{figure}

We have 
$$\bigcap_{P \in \mathcal{P}_{(u,w)}} \alpha(P) = \langle 1 \rangle ,\quad \bigcap_{P \in \mathcal{P}_{(u,z)}} \alpha(P) = \langle 2 \rangle, \;\text{ and}\:\: \bigcap_{P \in \mathcal{P}_{(w,z)}} \alpha(P) = \langle 3 \rangle,$$
and $\langle 2 \rangle + \langle 3 \rangle = \langle 1 \rangle$. Thus by Theorem \ref{pasting}, UDP holds for $(G, \alpha)$.
\end{example}

Next we will exhibit an edge labeled graph $(G, \alpha)$ that satisfies the hypotheses of Theorem \ref{pasting} but does not satisfy Equation \eqref{importantequation}.  This will be our first example of a graph on which the Universal Difference Property does not hold.  Before considering this example, we need to define a Prüfer domain.

\begin{definition}
\cite{Jensen} A \textit{Prüfer domain} is a commutative ring without zero divisors in which every non-zero finitely generated ideal is invertible. 
\end{definition}

\begin{proposition} \label{prop:5.10}\cite{Jensen}
Let $R$ be a domain. The following conditions are equivalent:
\begin{enumerate}
    \item $R$ is a Prüfer domain.
    \item If $I$, $J$, and $K$ are non-zero ideals of $R$, then $I + (J \cap K) = (I + J) \cap (I + K)$.
    \item If $I$, $J$, and $K$ are non-zero ideals of $R$, then $I \cap (J + K) = (I \cap J) + (I \cap K)$.
\end{enumerate}
\end{proposition}

We now use a base ring that is not a Prüfer domain to exhibit a family of pasted graphs on which the Universal Difference Property does not hold.

\begin{theorem}\label{Prufer domain counterexample}
Let $R$ be a ring which is not a Prüfer domain. Then there exists an edge labeled graph $(G, \alpha)$ over $R$ that does not satisfy UDP. 
\end{theorem}

\begin{figure}[H]
\begin{center}
\scalebox{0.2}{\includegraphics{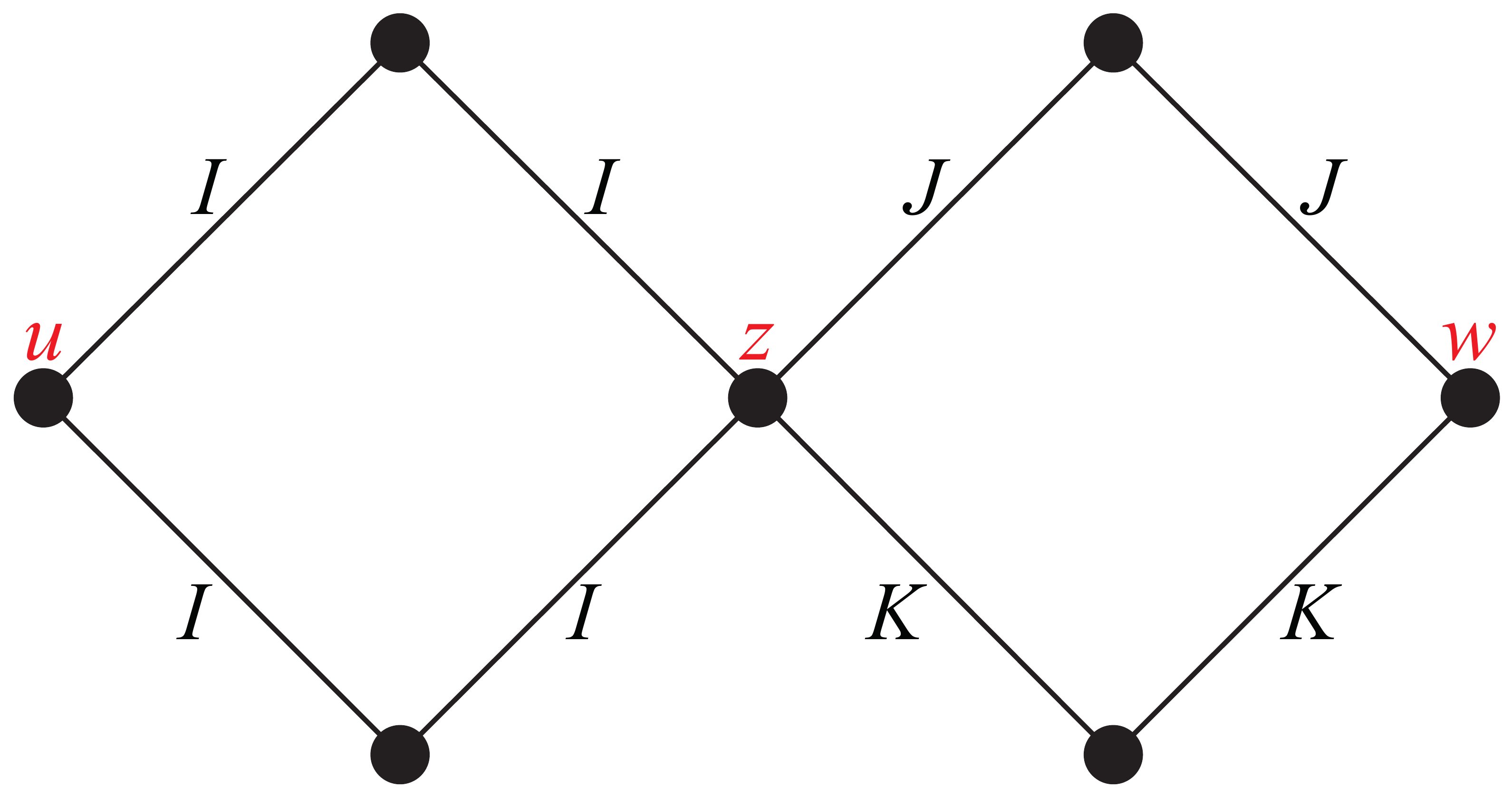}}
\caption{A graph labeled with certain ideals from a non-Prüfer domain}
\label{fig:generalnotprufer}
\end{center}
\end{figure}

\begin{proof}
Because $R$ is not a Prüfer domain, there exist non-zero ideals $I, J, K$ of $R$ such that $I + (J \cap K) \neq (I + J) \cap (I + K)$. Let $G_1$ and $G_2$ be $4$-cycles such that $V(G_1) \cap V(G_2) = \{z\}$, and let $G$ be the graph formed by pasting $G_1$ and $G_2$ at $z$. Fix $u \in V(G_1)$ and  $w \in V(G_2)$ such that neither $u$ nor $w$ is $z$ nor is adjacent to $z$, as in Figure \ref{fig:generalnotprufer}.  By Theorem \ref{pasting}, to show that our graph does not satisfy UDP, it suffices to show that Equation \eqref{importantequation} does not hold.

Now we introduce the edge labeling.  Label all edges of $G_1$ with the ideal $I$.  In $G_2$, there are two edge disjoint paths from $z$ to $w$.  Along one path, label each edge with the ideal $J$.  Along the other path, label each edge with the ideal $K$.  This labeling is depicted in Figure \ref{fig:generalnotprufer}.

Now \[ \bigcap_{P \in \mathcal{P}_{(u, w)}} \alpha(P) = (I + J) \cap (I + K), \]
while
\[ \bigcap_{P \in \mathcal{P}_{(u, z)}} \alpha(P) = I\]
and 
\[ \bigcap_{P \in \mathcal{P}_{(z, w)}} \alpha(P) = J \cap K. \]
Thus
\begin{equation*} \bigcap_{P \in \mathcal{P}_{(u, z)}} \alpha(P) + \bigcap_{P \in \mathcal{P}_{(z, w)}} \alpha(P) =I + (J \cap K) \neq(I + J) \cap (I + K)=\bigcap_{P \in \mathcal{P}_{(u, w)} } \alpha(P). \label{notequality} \end{equation*}

Hence Equation \eqref{importantequation} does not hold on this graph, and by Theorem $\ref{pasting}$, $(G, \alpha)$ does not satisfy UDP.  \end{proof}

\begin{example}
For a more specific example, let $R = \mathbb{Z}[x]$ and $(G,\alpha)$ be as in Figure \ref{fig:notprufer}. Note that $\mathbb{Z}[x]$ is not a Prüfer domain.  We have
\begin{equation*}
\bigcap_{P \in \mathcal{P}_{(u,w)}} \alpha(P) = \left(\langle 3 \rangle + \langle x + 3 \rangle\right) \cap \left(\langle 2 \rangle + \langle x - 3 \rangle\right) \cap \left(\langle 3 \rangle + \langle x - 3 \rangle\right) \cap \left(\langle 2 \rangle + \langle x + 3 \rangle\right),
\end{equation*}
while
\begin{equation*}
\bigcap_{P \in \mathcal{P}_{(u,z)}} \alpha(P) = \langle 3 \rangle \cap \langle 2 \rangle = \langle 6 \rangle
\end{equation*}
and
\begin{equation*}
\bigcap_{P \in \mathcal{P}_{(z,w)}} \alpha(P) = \langle x + 3 \rangle \cap \langle x - 3 \rangle = \langle x^2 - 9 \rangle.
\end{equation*}

 \noindent Note that $x + 3 \in \bigcap_{P \in \mathcal{P}_{(u,w)}} \alpha(P)$ but $x + 3 \notin \bigcap_{P \in \mathcal{P}_{(u,z)}} \alpha(P)$ and $x + 3 \notin \bigcap_{P \in \mathcal{P}_{(z,w)}} \alpha(P) $. Furthermore, $x + 3 \notin \langle 6 \rangle + \langle x^2 - 9 \rangle$, and so by Theorem \ref{pasting}, UDP does not hold on the edge labeled graph in Figure \ref{fig:notprufer}.

\begin{figure}[H]
\begin{center}
\scalebox{0.2}{\includegraphics{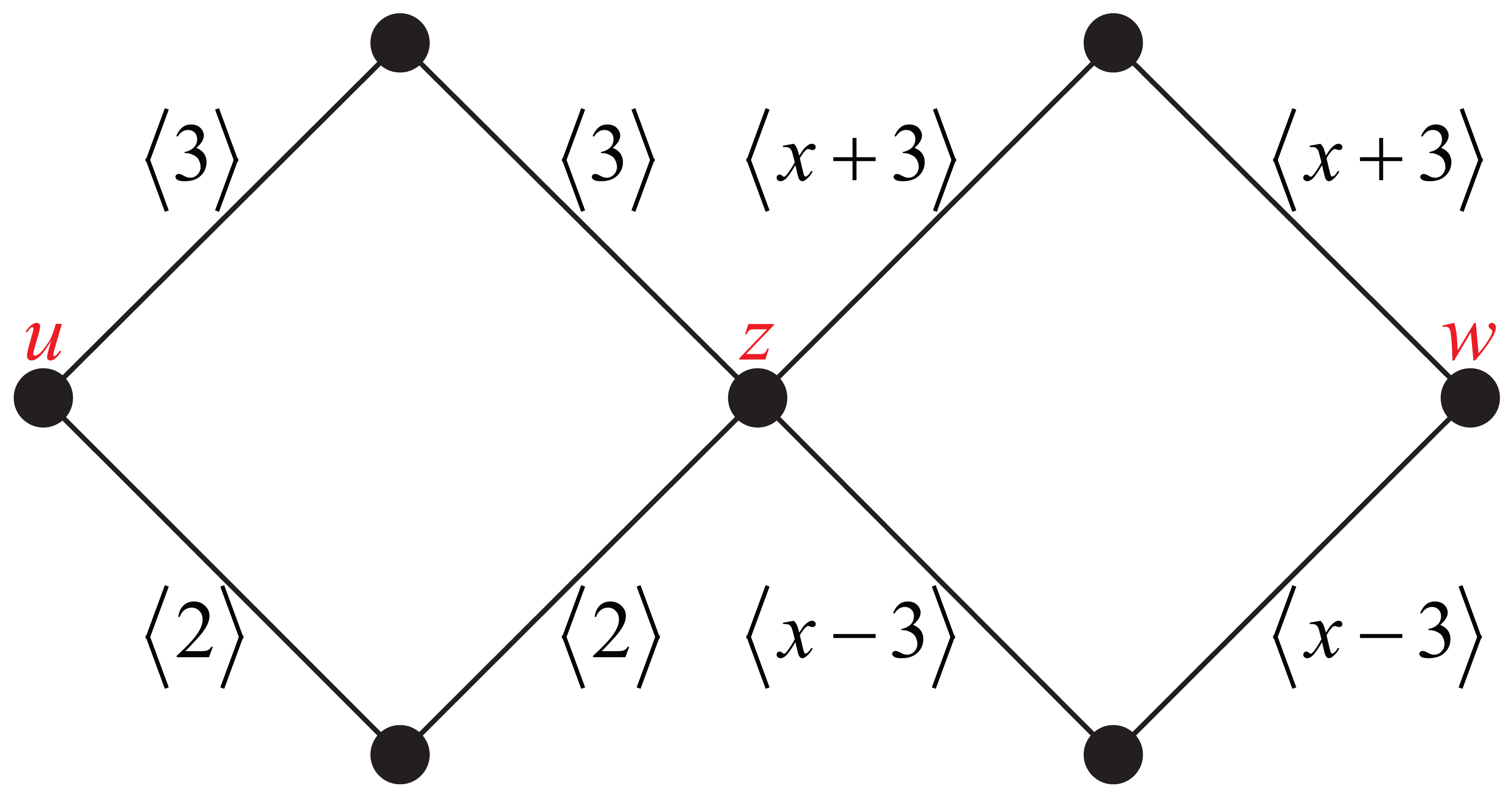}}
\caption{A graph with edge ideals of $\mathbb{Z}[x]$, a non-Prüfer domain}
\label{fig:notprufer}
\end{center}
\end{figure}

\end{example}

%%%%%%%%%%%%%%%%%%%%%%%%%%%%%%%%%%%
%%%% anticipated start of new section
%%%%%%%%%%%%%%%%%%%%%%%%%%%%%%%%%%%

\section{UDP over a Prüfer domain}\label{UDP over a Prufer domain}

The contrapositive of Theorem \ref{Prufer domain counterexample} is the statement, ``If $R$ is a ring such that every edge labeled graph over $R$ satisfies UDP, then $R$ is a Prüfer domain.''  In this section, we prove the converse of this contrapositive, to conclude that for a ring $R$, every edge labeled graph over $R$ satisfies UDP if and only if $R$ is a Prüfer domain.  We begin with a definition, which we will then use to give an alternate characterization of Prüfer domains.

\begin{definition}
The \textit{Elementwise Chinese Remainder Theorem} holds over a ring $R$ if for any elements $x_1,\ldots, x_n$ in $R$ and any ideals $I_1,\ldots, I_n$ of $R$,  the system of congruences
\begin{align*}
x&\equiv x_1 \bmod I_1\\
x&\equiv x_2 \bmod I_2\\
&\,\quad\vdots\\
x&\equiv x_n \bmod I_n\\
\end{align*}
has a solution in $R$ if and only if $x_j - x_k\in I_j+I_k$ for all $j,k$.
\end{definition}

\begin{theorem}[Proposition 25.1 of \cite{Gilmer}]\label{Gilmer characterization}
A ring $R$ is a Prüfer domain if and only if the Elementwise Chinese Remainder Theorem holds over $R$.
\end{theorem}

Next we prove a lemma similar to Lemma \ref{prop:subset} but with the assumption that $R$ is a Prüfer domain rather than the assumption that we are working with particular vertices in a pasted graph.

\begin{lemma}\label{referee's lemma}
Let $(G,\alpha)$ be an edge labeled graph over a Prüfer domain $R$.  For any $u,v,w\in V(G)$, we have
\begin{equation}\label{containment in referee's lemma}
\left(\bigcap_{P \in \mathcal{P}_{(v,w)}} \alpha(P)\right) \subseteq \left(\bigcap_{P \in \mathcal{P}_{(u,v)}} \alpha(P)\right)+\left(\bigcap_{P \in \mathcal{P}_{(u,w)}} \alpha(P)\right).
\end{equation}
\end{lemma}

\begin{proof}
Let $R$ be a Prüfer domain and $(G,\alpha)$ be an edge labeled graph over $R$.  Let $u,v,w\in V(G)$.  Let $P_1,\ldots,P_m$ be all the paths in $G$ from $u$ to $v$ and $Q_1,\ldots,Q_n$ be all the paths in $G$ from $u$ to $w$.  Then, because $R$ is a Prüfer domain,
\begin{align*}
\left(\bigcap_{P \in \mathcal{P}_{(u,v)}} \alpha(P)\right)+\left(\bigcap_{P \in \mathcal{P}_{(u,w)}} \alpha(P)\right)&=\left(\bigcap_{i=1}^m\alpha(P_i)\right)+\left(\bigcap_{j=1}^n \alpha(Q_j)\right)\\
&=\bigcap_{j=1}^n\left(\bigcap_{i=1}^m \alpha(P_i)+\alpha(Q_j)\right)\\
&=\bigcap_{j=1}^n \bigcap_{i=1}^m \left(\alpha(P_i)+\alpha(Q_j)\right).\\
\end{align*}

Now, to show the containment in \eqref{containment in referee's lemma}, we take an arbitrary element $x$ of $\bigcap_{P \in \mathcal{P}_{(v,w)}} \alpha(P)$ and will show that $x$ must be an element of $\bigcap_{j=1}^n \bigcap_{i=1}^m \left(\alpha(P_i)+\alpha(Q_j)\right)$.   Fix $1\leq j\leq n$ and $1\leq i\leq m$.  We will show that there exists a path $P$ from $v$ to $w$ such that $\alpha(P)\subseteq \alpha(P_i)+\alpha(Q_j)$.  Because $P_i$ is a path from $u$ to $v$ and $Q_j$ is a path from $u$ to $w$, concatenating $P_i$ and $Q_j$ yields a walk $W$ from $v$ to $w$, but this walk need not be a path from $v$ to $w$.  We can, however, take a subset of the edges and vertices in this walk $W$ to obtain a path $P$ from $v$ to $w$.  Since the edges in $P$ are a subset of the edges in $P_i$ and $Q_j$, we know $\alpha(P)\subseteq \alpha(P_i)+\alpha(Q_j)$.  Because $x$ belongs to $\alpha(P)$ for every path $P$ from $v$ to $w$, we know that $x$ belongs to $\alpha(P)$ for this particular path $P$ from $v$ to $w$, and thus $x$ must belong to $\alpha(P_i)+\alpha(Q_j)$.  Because $i$ and $j$ were arbitrary, we may conclude that $x$ belongs to the double intersection.
\end{proof}

We are now prepared to prove the main result of this section, a complete characterization of the set of rings over which every edge labeled graph must satisfy UDP.

\begin{theorem}\label{referee's big result}
For a ring $R$, every edge labeled graph over $R$ satisfies UDP if and only if $R$ is a Prüfer domain.
\end{theorem}

\begin{proof}
Let $R$ be a ring.  We have already shown in Theorem \ref{Prufer domain counterexample} that if every edge labeled graph over $R$ satisfies UDP, then $R$ is a Prüfer domain.  Here we show that if $R$ is a Prüfer domain, then every edge labeled graph over $R$ satisfies UDP.  Let $R$ be a Prüfer domain and $(G,\alpha)$ be an edge labeled graph over $R$.  Let $u,v\in V(G)$ and let $y\in\bigcap_{P \in \mathcal{P}_{(u,v)}} \alpha(P)$.  We will construct a spline $\rho:V(G)\rightarrow R$ such that $\rho(u)-\rho(v)=y$.

Label the vertices of $G$ as $v_0=u, v_1, v_2,\ldots,v_n=v$.  We begin defining the function $\rho$ from $V(G)$ to $R$ by letting $\rho(v_0)=y$ and $\rho(v_n)=0$.  Then $\rho(u)-\rho(v)=\rho(v_0)-\rho(v_n)=y$, as desired.  We will construct $\rho(v_1),\ldots,\rho(v_{n-1})$ by induction.

Our inductive hypothesis $H_\ell$ is that $\rho(v_j)-\rho(v_k)\in\bigcap_{P \in \mathcal{P}_{(v_j,v_k)}} \alpha(P)$ for all $0\leq j<k\leq \ell$ and $\rho(v_j)\equiv y\bmod \bigcap_{P \in \mathcal{P}_{(v_j,v)}} \alpha(P)$ for all $0\leq j\leq \ell$.  As our base case, notice that $H_0$ holds since $\rho(v_0)=y\equiv y\bmod \bigcap_{P \in \mathcal{P}_{(v_0,v)}} \alpha(P)$.  Let $1\leq i\leq n-1$ and assume $H_{i-1}$ holds.  We will show $H_i$ holds.  Because $H_{i-1}$ holds, we know that $\rho(v_j)-\rho(v_k)\in\bigcap_{P \in \mathcal{P}_{(v_j,v_k)}} \alpha(P)$ for all $0\leq j<k\leq i-1$ and $\rho(v_j)\equiv y\bmod \bigcap_{P \in \mathcal{P}_{(v_j,v)}} \alpha(P)$ for all $0\leq j\leq i-1$.  To conclude that $H_i$ holds, it remains to show that (i) $\rho(v_i)\equiv y\bmod \bigcap_{P \in \mathcal{P}_{(v_i,v)}} \alpha(P)$ and (ii) $\rho(v_j)-\rho(v_i)\in\bigcap_{P \in \mathcal{P}_{(v_j,v_i)}} \alpha(P)$ for all $0\leq j\leq i-1$. Item (ii) is the list of conditions
\begin{align*}
\rho(v_0)-\rho(v_i)&\in \bigcap_{P \in \mathcal{P}_{(v_0,v_i)}} \alpha(P),\\
\rho(v_1)-\rho(v_i)&\in \bigcap_{P \in \mathcal{P}_{(v_1,v_i)}} \alpha(P),\\
\rho(v_2)-\rho(v_i)&\in \bigcap_{P \in \mathcal{P}_{(v_2,v_i)}} \alpha(P),\\
&\;\;\vdots\\
\rho(v_{i-1})-\rho(v_i)&\in \bigcap_{P \in \mathcal{P}_{(v_{i-1},v_i)}} \alpha(P).\\
\end{align*}

Consider the following congruences:
\begin{align*}
x&\equiv \rho(v_0)\mod\bigcap_{P \in \mathcal{P}_{(v_0,v_i)}} \alpha(P),\\
x&\equiv\rho(v_1)\mod\bigcap_{P \in \mathcal{P}_{(v_1,v_i)}} \alpha(P),\\
x&\equiv\rho(v_2)\mod\bigcap_{P \in \mathcal{P}_{(v_2,v_i)}} \alpha(P),\\
&\;\;\vdots\\
x&\equiv\rho(v_{i-1})\mod\bigcap_{P \in \mathcal{P}_{(v_{i-1},v_i)}} \alpha(P).\\
\end{align*}
By our inductive hypothesis, 
\[
\rho(v_j)-\rho(v_k)\in\bigcap_{P \in \mathcal{P}_{(v_j,v_k)}} \alpha(P)\subseteq\left(\bigcap_{P \in \mathcal{P}_{(v_j,v_i)}} \alpha(P)\right)+\left(\bigcap_{P \in \mathcal{P}_{(v_k,v_i)}} \alpha(P)\right),
\]
for all $0\leq j<k\leq i-1$, where the last containment follows from Lemma \ref{referee's lemma}.  Also by the inductive hypothesis, 
\[
\rho(v_j)-y=\rho(v_j)-\rho(v_0)\in\bigcap_{P \in \mathcal{P}_{(v_j,v)}} \alpha(P)\subseteq\left(\bigcap_{P \in \mathcal{P}_{(v_i,v)}} \alpha(P)\right)+\left(\bigcap_{P \in \mathcal{P}_{(v_j,v_i)}} \alpha(P)\right),
\]
for all $0\leq j\leq i-1$.  Because $R$ is a Prüfer domain, we know by Theorem \ref{Gilmer characterization} that the Elementwise Chinese Remainder Theorem holds over $R$.  Combining this with the containments immediately above, we see that the system of congruences must have a solution in $R$.  We define $\rho(v_i)$ to be a solution to this system of congruences, and Items (i) and (ii) are both satisfied.  By induction, we can continue in this manner until we have defined $\rho(v_{n-1})$, and then we stop.

Now we have defined $\rho$ on all of $V(G)$, and we will show that $\rho$ is a spline.  Pick any pair of adjacent vertices $v_s$ and $v_t$ and assume $s<t$.  By construction, for any $0\leq j<k\leq n$, $\rho(v_k)-\rho(v_j)\in\bigcap_{P\in\mathcal{P}_{(v_k,v_j)}}\alpha(P)$.  Since one of the paths from $v_s$ to $v_t$ consists of the single edge $e_{st}$ from $v_s$ to $v_t$, it follows that $\bigcap_{P\in\mathcal{P}_{(v_s,v_t)}}\alpha(P)\subset \alpha(e_{st})$.  Thus $\rho(v_t)-\rho(v_s)\in\alpha(e_{st})$, so $\rho$ is a spline on $(G,\alpha)$, and we already knew that $\rho(u)-\rho(v)=y$.
\end{proof}

A Bezout domain is an integral domain such that any finitely generated ideal is principal. Every PID is a Bezout domain, and every Bezout domain is a Prüfer domain.  Consequently, by Theorem \ref{referee's big result}, we know that any edge labeled graph over a Bezout domain satisfies UDP.  In particular, any edge labeled graph over a PID satistifes UDP.  
\section{UDP is a structural property of an edge labeled graph}
\label{UDP and isomorphisms}

In this section we show that the Universal Difference Property is a structural property of an edge labeled graph; that is, the Universal Difference Property is preserved by an isomorphism of edge labeled graphs.  We begin with a definition.
\begin{definition}\label{def1}
\cite{Gilbert} Let $(G,\alpha)$ and $(G',\alpha')$ be edge labeled graphs with base ring $R$. A \textit{homomorphism of edge labeled graphs} $\varphi: (G,\alpha) \to (G',\alpha')$ is a graph homomorphism $\varphi_1: G \to G'$ together with a ring automorphism $\varphi_2:R \to R$ so that for each edge $e \in E(G)$ we have $\varphi_2(\alpha(e)) = \alpha'(\varphi_1(e))$. An \textit{isomorphism of edge labeled graphs} is a homomorphism of edge labeled graphs whose underlying graph homomorphism is in fact an isomorphism.
\end{definition}

We are now ready to state and prove the main result of this section.

\begin{theorem}\label{SIUDP}
Let $(G,\alpha)$ and $(G',\alpha')$ be edge labeled isomorphic. If the Universal Difference Property holds on $(G,\alpha)$, then it holds on $(G^{'},\alpha')$.
\end{theorem}
\begin{proof}
Let $u',w'\in V(G')$, and $P_{1}',P_{2}',\ldots,P_{n}'$ be the paths in $G'$ from $u'$ to $w'$. Let $y\in\bigcap_{i=1}^{n}\alpha'(P_{i}')$. Let $\varphi_{1}:G\rightarrow G'$ be the graph isomorphism and $\varphi_{2}:R\rightarrow R$ be the ring automorphism from Definition \ref{def1}. Let $u,v\in V(G)$ such that $\varphi_{1}^{-1}(u')=u$ and $\varphi_{1}^{-1}(w')=w$. Finally, let $\alpha$ and $\alpha'$ be the edge labelings of $G$ and $G'$, respectively. We know there must be a $z\in R$ such that $\varphi_{2}(z)=y$.

Consider the path $P_{i}'$ from $u'$ to $w'$. Let $e_{1}',e_{2}',\ldots,e_{k}'$ be the edges in the path $P_{i}'$. We can write $y$ as a sum of elements from the ideals along $P_{i}'$ such that $y=a_{1}'+a_{2}'+\cdots+a_{k}'$, where $a_{i}'\in\alpha'(e_{i}')$.\\
Then we have
\begin{align*}
       z&=\varphi_{2}^{-1}(y)\\
        &=\varphi_{2}^{-1}(a_{1}'+a_{2}'+\dots+a_{k}')\\
        &=\varphi_{2}^{-1}(a_{1}')+\varphi_{2}^{-1}(a_{2}')+\dots+\varphi_{2}^{-1}(a_{k}').
\end{align*}

Consider the element $a_{j}' \in \alpha '(e_j ')$.
We know there exists an edge $e_j$ in $G$ such that $\varphi_1 (e_j) = e_j '$.
Then $a_j ' \in \alpha '(e_j ') = \alpha '(\varphi_1 (e_j)) = \varphi_2 (\alpha(e_j))$ by Definition \ref{def1}.
Then $\varphi_{2}^{-1} (a_j ') \in \alpha(e_j)$.
Let  $\varphi_{2}^{-1} (a_j ') = a_j$.  Doing this for $1\leq j\leq k$, we obtain
\begin{align*}
       z&=\varphi_{2}^{-1}(y)\\
        &=\varphi_{2}^{-1}(a_{1}'+a_{2}'+\dots+a_{k}')\\
        &=\varphi_{2}^{-1}(a_{1}')+\varphi_{2}^{-1}(a_{2}')+\dots+\varphi_{2}^{-1}(a_{k}')\\ 
        &= a_1 + a_2 + \dots + a_k.
\end{align*}

Recalling that $e_1, e_2, \dots, e_k$ is a path in $G$ from $u$ to $w$, we have $z$ as a sum of elements of the ideals labeling the edges along this path.  Since vertex adjacency is preserved by the graph isomorphism, every path from $u$ to $w$ in $G$ corresponds to one of the paths $P_i'$ from $u'$ to $w'$ in $G'$. Hence  $z\in\bigcap_{P \in \mathcal{P}_{(u, w)}} \alpha(P)$.  Since $G$ satisfies the Universal Difference Property, there is a spline $\rho\in R_G$ such that $\rho(u)-\rho(w)=z$. Then \begin{align}\label{proofend} y =\varphi_{2}(z)=\varphi_{2}(\rho(u)-\rho(w))=\varphi_{2}(\rho(u))-\varphi_{2}(\rho(w)).\end{align}

For any $v' \in V(G')$, let $v$ be the element of $G$ such that $\varphi_1(v) = v'$. Define a vertex labeling $\gamma$ on $G'$ by $\gamma(v') = \varphi_2(\rho(\varphi_1^{-1}(v'))) = \varphi_2(\rho(v))$. We will show that $\gamma$ is a spline on $G'$ and that $\gamma(u') - \gamma(w') = y$. Let $s'$ and $t'$ be adjacent vertices in $G'$. Then there exist vertices $s$ and $t$ in $V(G)$ such that $\varphi_1(s) = s'$ and $\varphi_1(t) = t'$. Then
\begin{align*}
    \gamma(s') - \gamma(t') &= \varphi_2(\rho(\varphi_1^{-1}(s'))) - \varphi_2(\rho(\varphi_1^{-1}(t'))) \\
                            &= \varphi_2(\rho(s)) - \varphi_2(\rho(t)) \\
                            &= \varphi_2(\rho(s) - \rho(t)) \\
                            &\in \varphi_2(\alpha(st)),
\end{align*}
but $\varphi_2(\alpha(st)) = \alpha'(\varphi_1(st)) = \alpha'(s't')$ and so $\gamma$ is indeed a spline on $G'$.  Lastly, we have 
\begin{align*}
    \gamma(u')-\gamma(w')&=\varphi_{2}(\rho(u))-\varphi_{2}(\rho(w))\\
    &=y\quad \text{by Equation \eqref{proofend}}.
\end{align*}
\noindent Thus $(G',\alpha')$ satisfies UDP.
\end{proof}

\section{Acknowledgments}
The authors wish to thank Professor Julianna Tymoczko for her thoughtful feedback and questions.  The authors also wish to thank the referees for their time and helpful comments and in particular for showing us how to strengthen our results and obtain Theorem \ref{referee's big result} and its proof.

The second author thanks the Institute for Computational and Experimental Research in Mathematics and the American Institute of Mathematics for their generous support of the Research Experience for Undergraduate Faculty (REUF) program, in which she participated in Summers 2017 and 2018.  The second, third, fourth, and fifth authors acknowledge support from the University of Texas at Tyler Research Experience for Undergraduates,  NSF Grant DMS-1659221, in which the third, fourth, and fifth authors participated during Summer 2019, under the mentorship of the second author.  The second and seventh authors acknowledge support from NSF Grant HRD-1826745 Louis Stokes STEM Pathways and Research Alliance: University of Texas System LSAMP, in which the seventh author participated as an undergraduate researcher under the mentorship of the second author in the Summer Research Academy 2019.

The first and sixth authors are supported by Hacettepe University Scientific Research Projects Coordination Unit. Project Number: FHD-2022-19802.





\end{document}